# Universal derivative-free optimization method with quadratic convergence

Sergey N. Moiseev

Kodofon,

394062, Russia, Voronezh

smoiseev@kodofon.vrn.ru

2011

**Abstract**

A new universal derivative-free optimization method CDOS (Conjugate Direction with Orthogonal Shift) is proposed. The CDOS method was specially developed to solve optimization tasks where the objective function and constraints are black boxes. The method has quadratic convergence for quadratic and near quadratic functions. An objective function can be non-differentiable and non-continuous. For constrained optimization the constraints can also be non-differentiable and non-continuous. The method handles inequality constraints directly, i.e., it does not transform the objective function, nor uses numeric values of the inequality constraints – it uses only the fact of constraint violation.

**Keywords:** non-linear optimization, derivative-free methods, conjugate directions, quadratic convergence, constrained optimization, global optimization, multimodal optimization, mixed integer-discrete-continuous optimization

# Introduction

In practice many objective functions are non-differentiable or even non-continuous. There are objective functions in the form of complicated programs, i.e., in the form of black boxes. These objective functions can even be differentiable but you cannot get the derivatives of these functions. Fast optimization methods using derivatives cannot be used for such objective functions. In these cases direct search derivative-free methods are usually used, for example, the Nelder-Mead simplex method. Though the Nelder-Mead method is good enough for non-differentiable functions with small number of variables, it does not have quadratic convergence for quadratic functions. Moreover, the Nelder-Mead method is too slow and unreliable for high dimension optimization problems.

There are derivative-free optimization methods having quadratic convergence for quadratic and near quadratic functions such as Powell's and Brent's line-search methods. However, these methods are highly unreliable for non-differentiable, non-continuous functions and for constrained optimization tasks.

The optimization constraints can also be non-differentiable and non-continuous. In some cases the constraints may not even be given explicitly. For example, let an objective function be a complex program simulating the telecommunication network. Such simulator has many parameters. The constraints for these parameters are included into the simulator by a programmer and are not available explicitly. You can only get some error message when the constraints are violated. Fast regular optimization methods using derivatives, for example, SQP method (Sequential Quadratic Programming), cannot be used for these constrained optimization tasks.

In this paper we present a new universal derivative-free line-search CDOS (Conjugate Direction with Orthogonal Shift) optimization method. The CDOS method was specially developed to solve the optimization tasks where the objective function and inequality constraints are black boxes and the objective function is not defined outside inequality constraint boundaries. This method is good for all the above mentioned optimization problems. The CDOS method is also fast for differentiable functions because it has quadratic convergence for quadratic and near quadratic functions. The method requires minimal possible $\frac{n(n+1)}{2}$ one-dimensional line searches to calculate exactly the optima of $n$-dimensional quadratic function.

The CDOS method is not quite greedy, that is, moving in directions that do not contain the current extremum point is allowed. Therefore, it can be used as a basis for creating effective global optimization methods.

The CDOS method has the following features:

- it has quadratic convergence for quadratic and near quadratic functions,
- it is more reliable than the Nelder-Mead method for non-differentiable functions,
- it is able to solve the high dimensional optimization problem,
- it does not require the objective function to be defined outside the constraint boundaries,
- it can handle inequality constraints directly,
- it moves fast along the boundaries of the constraints,
- it can be easily transformed into a global optima searcher because it is not a greedy searcher,
- it can be easily transformed into a multimodal optimizer because it has quadratic convergence,
- it can be easily transformed into a mixed integer-discrete-continuous optimizer because it is not a greedy searcher and has quadratic convergence,
- it is easily programmed.

## Optimization problem formulation

The optimization problem can be formulated as follows:

minimize (maximize)

$$f(\mathbf{x}), \mathbf{x} \in R^n,$$

subject to

$$\mathbf{c}_0(\mathbf{x}) \rightarrow \{\text{true, false}\},$$
$$\mathbf{c}_1(\mathbf{x}) \geq 0,$$
$$\mathbf{c}_2(\mathbf{x}) > 0,$$
$$\mathbf{c}_3(\mathbf{x}) \neq 0,$$
$$\mathbf{h}(\mathbf{x}) = 0,$$

where $\mathbf{x}$ are the optimization problem variables, $f$ is the objective function, $\mathbf{c}_i$ are the inequality constraint functions, $\mathbf{h}$ are the equality constraint functions. The constraints $\mathbf{c}_0$ describe only those constraints that were violated (false), or not (true). Note that inequality constraints $\mathbf{c}_0$ can be transformed into $\mathbf{c}_2$ type of constraints $\mathbf{c}_{02} > 0$ as follows: $\mathbf{c}_{02} = 1$ when $\mathbf{c}_0 = \text{true}$ and $\mathbf{c}_{02} = -1$ when $\mathbf{c}_0 = \text{false}$. The values $\{1, -1\}$ of inequality constraint functions $\mathbf{c}_{02}$ do not contain any information about the degree of constraints $\mathbf{c}_0$ violation.

# Method description

The CDOS method, as well as Powell's and Brent's methods, uses conjugate line-search directions. The two directions $\mathbf{u}_i$ and $\mathbf{u}_j$ are conjugated if $\mathbf{u}_i^T \mathbf{H} \mathbf{u}_j = 0$, $i \neq j$; $\mathbf{u}_i^T \mathbf{H} \mathbf{u}_j \geq 0$, $i = j$ where $\mathbf{H}$ is a positive definite Hessian matrix: $\mathbf{H} = \left| \frac{\partial}{\partial x_i} \left( \frac{\partial}{\partial x_j} f(\mathbf{x}) \right) \right|$, $i, j = 1..n$. If the function $f(x_1,...,x_n)$ is the $n$-dimensional quadratic function its minimum (maximum) is obtained for $n$ moves in $n$ various conjugated directions. Therefore the extremum point of the $n$-dimensional quadratic function is reached when $n$ conjugate search directions are constructed for the first time.

To construct the conjugate directions an orthogonal shift from the subset of already constructed conjugate directions is used. The orthogonal shift makes the algorithm non greedy but greatly increases its reliability and gives it new unique features. These new features lead to principal differences between the CDOS method and other methods with quadratic convergence. The most important new features are the fast moving along constraint boundaries and the drastically increasing reliability of non-differentiable function optimization.

The minimization algorithm consists of three stages. At stage one the first search direction is defined as being the opposite of quasi-gradient. At stage two the first $n$ conjugate directions using the orthogonal shift are constructed. At the third stage the main iteration cycle is carried out. This cycle covers the iterative update of the constructed conjugate directions using the orthogonal shift. Below we consider these three stages of minimum search for the $n$-dimensional function $f(x_1,...,x_n)$, $n > 1$ in more detail.

**Stage I.** Initial $n$ search directions correspond to coordinate axes: $\mathbf{u}_1 = (1, 0, 0,..., 0)^T$, $\mathbf{u}_2 = (0, 1, 0,..., 0)^T$, ..., $\mathbf{u}_n = (0, 0, 0,..., 1)^T$. One step is made in each direction starting from the initial point. The default value for the initial point is $\mathbf{x}_0 = (0.9,...,0.9)$. The initial step value $\lambda$ is determined by a user. By default it is equal to 1. The vector of the objective function increments for each direction $\mathbf{s} = (\Delta f_1,...,\Delta f_n)^T$ is calculated. As the first conjugate direction we select the anti-gradient direction which corresponds to the normalized vector $\mathbf{u}_1^* = -\frac{\mathbf{s}}{\|\mathbf{s}\|}$. The first search direction is replaced by the anti-gradient direction $\mathbf{u}_1 = \mathbf{u}_1^*$ and we carry out the one-dimension minimum search for this direction. Let $\mathbf{x}_{\min}^{(1)}$ denote the point of obtained minimum.

**Stage II.** Initial $n$ search directions at this stage are: $\mathbf{u}_1 = \mathbf{u}_1^*$, $\mathbf{u}_2 = (0, 1, 0,..., 0)^T$, …, $\mathbf{u}_n = (0, 0, 0,..., 1)^T$. The step value $\lambda$ at this stage is the same and equal to the initial step value. The orthogonal shift value is $\lambda_s = 0.62\lambda$ and does not change for the whole stage II. We shall describe constructing the rest of $n-1$ conjugate directions as the following pseudo-code.

**for i from 2 to n do**   # the cycle of creating 2, 3,…, $n$ conjugate directions

We make the orthogonal shift from the already available $\mathbf{u}_1, ..., \mathbf{u}_{i-1}$ conjugate directions from the current minimum point $\mathbf{x}_{\min}^{(i-1)}$. The orthogonal shift direction can be obtained by means of the Gram-Schmidt orthogonalization process or by QR decomposition of matrix $\mathbf{U} = (\mathbf{u}_1, ..., \mathbf{u}_i) = \mathbf{QR}$.

The Gram–Schmidt process takes a set of vectors $\mathbf{u}_1, ..., \mathbf{u}_i$ and generates an orthonormal set of vectors $\mathbf{u}_1^*, ..., \mathbf{u}_i^*$. Vectors $\mathbf{u}_1^*, ..., \mathbf{u}_{i-1}^*$ are not used. The orthogonal shift direction is defined by vector $\mathbf{u}_i^*$. This way we have a point corresponding to the orthogonal shift $\mathbf{y}_{\min}^{(i-1)} = \mathbf{x}_{\min}^{(i-1)} + \lambda_s \mathbf{u}_i^*$.

If QR decomposition is used then the columns of matrix $\mathbf{Q} = (\mathbf{q}_1,...,\mathbf{q}_i)$ are an orthonormal set of vectors and the orthogonal shift direction is defined by the last $i$-th column $\mathbf{q}_i$ of matrix $\mathbf{Q}$. Thus we have a point corresponding to the orthogonal shift $\mathbf{y}_{\min}^{(i-1)} = \mathbf{x}_{\min}^{(i-1)} + \lambda_s \mathbf{q}_i$.

The QR decomposition that uses Householder transformations is faster than the Gram-Schmidt process and is more preferable.

   **for j from 1 to i-1 do**   # repeated movement in $i-1$ conjugate directions

   Starting from point $\mathbf{y}_{\min}^{(i-1)}$ we perform a one-dimensional search for minimum in $\mathbf{u}_j$ direction. If the function value in the obtained minimum point is less than the function value in point $\mathbf{y}_{\min}^{(i-1)}$, we update point $\mathbf{y}_{\min}^{(i-1)}$.

   **end for j**.

   Direction $\mathbf{u}_i$ is replaced by a new direction connecting points $\mathbf{x}_{\min}^{(i-1)}$ and $\mathbf{y}_{\min}^{(i-1)}$ as follows:

$$\mathbf{u}_i = \begin{cases} \dfrac{\mathbf{x}_{\min}^{(i-1)} - \mathbf{y}_{\min}^{(i-1)}}{\|\mathbf{x}_{\min}^{(i-1)} - \mathbf{y}_{\min}^{(i-1)}\|}, & f(\mathbf{x}_{\min}^{(i-1)}) < f(\mathbf{y}_{\min}^{(i-1)}), \\ \dfrac{\mathbf{y}_{\min}^{(i-1)} - \mathbf{x}_{\min}^{(i-1)}}{\|\mathbf{x}_{\min}^{(i-1)} - \mathbf{y}_{\min}^{(i-1)}\|}, & f(\mathbf{x}_{\min}^{(i-1)}) \geq f(\mathbf{y}_{\min}^{(i-1)}). \end{cases}$$

This direction will be conjugate to all $\mathbf{u}_1, ..., \mathbf{u}_{i-1}$ previous conjugate directions. A one-dimensional search for minimum in $\mathbf{u}_i$ direction is performed starting from point $\mathbf{x}_{\min}^{(i-1)}$, if $f(\mathbf{x}_{\min}^{(i-1)}) < f(\mathbf{y}_{\min}^{(i-1)})$, or from point $\mathbf{y}_{\min}^{(i-1)}$ otherwise. The resulting minimum point becomes the current minimum point $\mathbf{x}_{\min}^{(i)}$ for the following iteration cycle.

**end for i**.

If the target function is the $n$-dimensional quadratic function, its minimum is reached already at the end of stage II.

**Stage III.** Initial $n$ search directions at this stage correspond to mutual conjugate directions found at stages II: $\mathbf{u}_1, ..., \mathbf{u}_n$. The initial step value $\lambda^{(1)} = 0.3 \|\mathbf{x}_{\min}^{(n)} - \mathbf{x}_{\min}^{(n-1)}\| + 0.091 \cdot \lambda$. If $\lambda^{(1)} = 0$, then $\lambda^{(1)} = \text{tol}$, where tol is a confidence interval for the extremum point. The initial orthogonal shift value is $\lambda_s^{(1)} = 0.62 \lambda^{(1)}$. If $\lambda_s^{(1)} = 0$, then $\lambda_s^{(1)} = \lambda^{(1)}$. Let $\mathbf{x}_{\min}^{(1)}$ (and $\mathbf{x}_{\min}^{(0)}$) denote the current minimum point. We shall describe the main iteration cycle as the following pseudo-code.

**for i from 1 to infinity do**     # main iteration cycle

We perform the orthogonal shift from the already available $\mathbf{u}_n, ..., \mathbf{u}_2$ conjugate directions from the current minimum point $\mathbf{x}_{\min}^{(i)}$. The orthogonal shift direction can be obtained by means of the Gram-Schmidt orthogonalization process or by QR decomposition of matrix $\mathbf{U} = (\mathbf{u}_n, ..., \mathbf{u}_1) = \mathbf{QR}$.

The Gram–Schmidt process takes a set of vectors $\mathbf{u}_n, ..., \mathbf{u}_1$ and generates an orthonormal set of vectors $\mathbf{u}_n^*, ..., \mathbf{u}_1^*$. Vectors $\mathbf{u}_n^*, ..., \mathbf{u}_2^*$ are not used. The orthogonal shift direction is defined by vector $\mathbf{u}_1^*$. This way we have a point corresponding to the orthogonal shift $\mathbf{y}_{\min}^{(i)} = \mathbf{x}_{\min}^{(i)} + \lambda_s^{(i)} \mathbf{u}_1^*$.

If QR decomposition is used then $\mathbf{Q} = (\mathbf{q}_1, ..., \mathbf{q}_n)$ is an orthogonal matrix (its columns are orthogonal unit vectors meaning $\mathbf{Q}^T \mathbf{Q} = \mathbf{I}$) and the orthogonal shift direction is defined by the last column $\mathbf{q}_n$ of matrix $\mathbf{Q}$. If $\|\mathbf{q}_n - \mathbf{u}_1\| > \|\mathbf{q}_n + \mathbf{u}_1\|$ we change the sign of the direction $\mathbf{q}_n$: $\mathbf{q}_n = -\mathbf{q}_n$. Thus we have a point corresponding to the orthogonal shift $\mathbf{y}_{\min}^{(i-1)} = \mathbf{x}_{\min}^{(i-1)} + \lambda_s \mathbf{q}_n$.

The QR decomposition using Householder transformations is faster than the Gram-Schmidt process and is more preferable.

We circularly shift the whole set of conjugate directions $\mathbf{u}_1, ..., \mathbf{u}_n$ to the left. In other words we re-denote the set of conjugate directions as follows: let us denote $\mathbf{u}_2$ as $\mathbf{u}_1$, $\mathbf{u}_3$ as $\mathbf{u}_2$, ..., $\mathbf{u}_n$ as $\mathbf{u}_{n-1}$ and $\mathbf{u}_1$ as $\mathbf{u}_n$.

**for j from 1 to n-1 do**   # repeated movement in $n-1$ conjugate directions

    The one-dimensional search for minimum within this cycle shall be performed with a $3\lambda^{(i)}$ step size. This is done because we have deviated from the current minimum and perform a non-greedy search. Therefore the typical scale of function variation has increased.

    Starting from point $\mathbf{y}_{\min}^{(i)}$ we perform a one-dimensional search for minimum in $\mathbf{u}_j$ direction. If the function value in the obtained point is less than the function value in point $\mathbf{y}_{\min}^{(i)}$, we update point $\mathbf{y}_{\min}^{(i)}$.

**end for j**.

Direction $\mathbf{u}_n$ is replaced by a new direction connecting points $\mathbf{x}_{\min}^{(i)}$ and $\mathbf{y}_{\min}^{(i)}$ as follows:

$$\mathbf{u}_n = \begin{cases} \dfrac{\mathbf{x}_{\min}^{(i)} - \mathbf{y}_{\min}^{(i)}}{\|\mathbf{x}_{\min}^{(i)} - \mathbf{y}_{\min}^{(i)}\|}, & f(\mathbf{x}_{\min}^{(i)}) < f(\mathbf{y}_{\min}^{(i)}), \\ \dfrac{\mathbf{y}_{\min}^{(i)} - \mathbf{x}_{\min}^{(i)}}{\|\mathbf{x}_{\min}^{(i)} - \mathbf{y}_{\min}^{(i)}\|}, & f(\mathbf{x}_{\min}^{(i)}) \geq f(\mathbf{y}_{\min}^{(i)}). \end{cases}$$

This direction will be conjugate to all $\mathbf{u}_1, ..., \mathbf{u}_{n-1}$ previous conjugate directions. A one-dimensional search for minimum in $\mathbf{u}_n$ direction is performed with step $\lambda^{(i)}$ starting from point $\mathbf{x}_{\min}^{(i)}$, if $f(\mathbf{x}_{\min}^{(i)}) < f(\mathbf{y}_{\min}^{(i)})$, or from point $\mathbf{y}_{\min}^{(i)}$ otherwise. If the function value in the obtained point is less than the function value in point $\mathbf{x}_{\min}^{(i)}$, we update point $\mathbf{x}_{\min}^{(i)}$.

The search step and orthogonal shift values are adjusted for the following iteration as follows: $\lambda^{(i+1)} = 0.3 \|\mathbf{x}_{\min}^{(i)} - \mathbf{x}_{\min}^{(i-1)}\| + 0.091 \cdot \lambda^{(i)}$. If $\lambda^{(i+1)} = 0$, then $\lambda^{(i+1)} = \text{tol}$, where tol is a confidence interval for the extremum point. The orthogonal shift value is $\lambda_s^{(i+1)} = 0.62 \lambda^{(i+1)}$. If $\lambda_s^{(i+1)} = 0$, then $\lambda_s^{(i+1)} = \lambda^{(i+1)}$.

We shall now check the following condition of search termination. If the following inequalities are satisfied $N_{exit}$ times in a row: $\lambda^{(i+1)} \leq tol$ and $f(\mathbf{x}_{min}^{(i-1)}) - f(\mathbf{x}_{min}^{(i)}) \leq tol_2$, the search is stopped. Here $tol_2$ is a confidence interval value for the extremum value. For the majority function $N_{exit} = 2$ is enough. For non-differentiable functions the value $N_{exit} = 10$ is more reliable. As a minimum estimate we use $f(\mathbf{x}_{min}^{(i)})$, and as a minimum point we use $\mathbf{x}_{min}^{(i)}$. Otherwise we move on to the next search iteration. The resulting minimum point becomes the current minimum point $\mathbf{x}_{min}^{(i+1)}$ for the following iteration cycle.

**end for i**.

Fig. 1 below shows a real path of the following function minimum search
$$f(x, y) = x^2 + y^2 - 1.5xy,$$
starting with initial point $[x = 5, y = 3]$. The real search is performed by Search command from DirectSearch package, version 2 [1]. The circles denote the trial points.

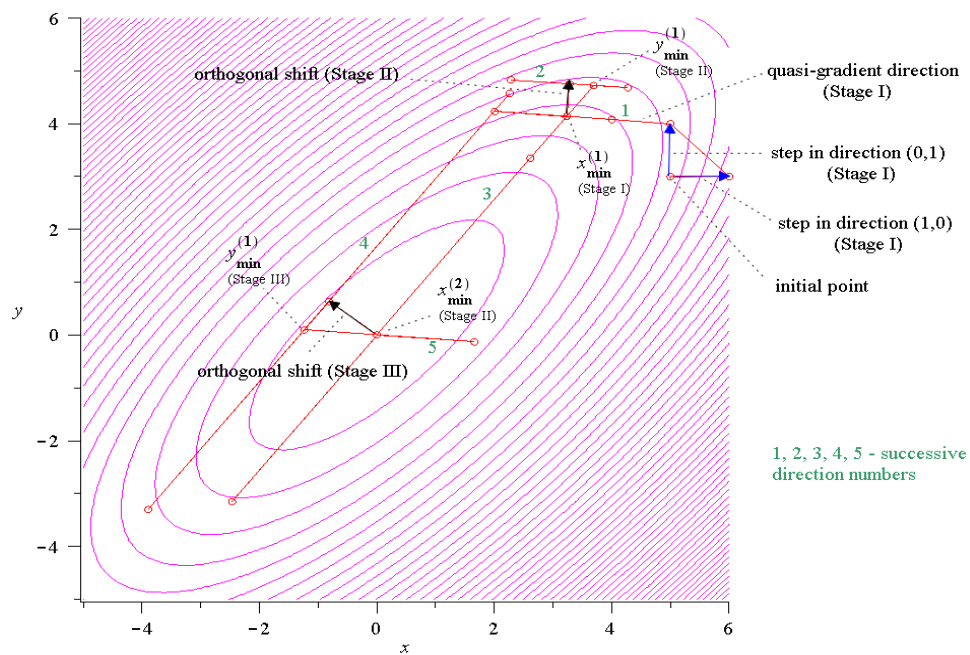

**Fig. 1. Search path**

The above plot shows that the quadratic function minimum is reached by the end of Stage II by three line-direction searches. One iteration of Stage III leads to the same minimum point obtained at Stage II.

# One dimensional line-search algorithm

A one-dimensional search in case $n > 1$ is carried out as follows. A step is made from the initial point in the search direction. If the step is successful, the extremum point is updated, the step is doubled and the search is continued in the same direction until we have the first unsuccessful step. Then a single parabolic approximation is performed for the last three points. As a minimum point we take a point where the minimum function value is achieved. If the first step is unsuccessful, we change the search direction to the opposite direction and repeat this procedure.

The parabolic approximation is performed as follows. Let $\mathbf{x}', \mathbf{x}'', \mathbf{x}'''$ be the three distinguishing points in one line-direction $\mathbf{u}_1$. Let $d_{13} > d_{12}, d_{13} > d_{23}$, where $d_{13} = \|\mathbf{x}' - \mathbf{x}'''\|, d_{12} = \|\mathbf{x}' - \mathbf{x}''\|, d_{23} = \|\mathbf{x}'' - \mathbf{x}'''\|$. Then the minimum point of the parabolic approximation is

$$\mathbf{x}_{\min} = \mathbf{x}' + d \cdot \mathbf{u}_1,$$

where

$$d = \frac{1}{2} \frac{x_0^2 (f''' - f'') + x_1^2 (f' - f''') + x_2^2 (f'' - f')}{f'''(x_0 - x_1) + f''(x_2 - x_0) + f'(x_1 - x_2)},$$
$$f' = f(\mathbf{x}'), f'' = f(\mathbf{x}''), f''' = f(\mathbf{x}'''),$$
$$x_0 = 0, x_1 = d_{12}, x_2 = d_{13}.$$

# Search in space-curve direction

In order to effectively move along the curved valleys, such as in the Rosenbrock function, CDOS uses the search in space-curve directions. One method describing how to use the search in a curved direction is described in [4]. The CDOS uses another method as follows.

Let $\mathbf{x}' = (x'_1, \ldots, x'_n), \mathbf{x}'' = (x''_1, \ldots, x''_n), \mathbf{x}''' = (x'''_1, \ldots, x'''_n)$ be the three successive distinguishing approximations of the minimum point that were obtained during the $i$-th, $(i + n + 1)$-th and $(i + 2n + 2)$-th iterations. It is clear that $f(\mathbf{x}''') \leq f(\mathbf{x}'') \leq f(\mathbf{x}')$. From all $n$ coordinates we select such $m$-th coordinate for which we have

$$x'_m < x''_m < x'''_m \text{ or } x'_m > x''_m > x'''_m.$$

Now we consider $x'_m, x''_m, x'''_m$ as argument values of the parabolic function and $x'_k, x''_k, x'''_k, k \neq m$ as the parabolic function values that correspond to these argument values. For each $k \neq m$ a parabola is uniquely defined by these three points. So we use $n - 1$ parabolic extrapolations to move in a curved direction. For example, when $x'_m < x''_m < x'''_m$ the new curved direction searching point is

$$\mathbf{x}_c = (x_1, \ldots, x_{m-1}, x'''_m + \lambda_c, x_{m+1}, \ldots, x_n),$$

where $x_k, k \neq m$ are corresponding parabola values for argument $x'''_m + \lambda_c$, $\lambda_c > 0$ is the curved direction step size. When $x'_m > x''_m > x'''_m$ the new curved direction searching point is

$$\mathbf{x}_c = (x_1, \ldots, x_{m-1}, x'''_m - \lambda_c, x_{m+1}, \ldots, x_n).$$

If the function value $f(\mathbf{x}_c)$ in the obtained point is less than the function value in the current minimum point $\mathbf{x}_{\min}^{(i)}$, we update point $\mathbf{x}_{\min}^{(i)}$.

## Constrained optimization

The CDOS method handles inequality constraints directly, that is, it does not transform an objective function, nor it uses numeric values of inequality constraints – it uses only the fact of constraint violation. So the CDOS method never computes the function values that do not satisfy inequality constraints, instead it searches for the feasible point.

Let us perform the minimum search in $\mathbf{u}_k$ direction starting from feasible point $\mathbf{x}$ with step $\lambda_0$. Assume that point $\mathbf{x} + \lambda_0 \mathbf{u}_k$ is an infeasible point. In this case we reduce the step size to $\lambda_1 = \varepsilon_1 \lambda_0$, $0 < \varepsilon_1 < 1$ and try another point $\mathbf{x} + \lambda_1 \mathbf{u}_k$. If this point is also infeasible we repeat the process. If we do not find a feasible point after some trial we change the search direction to the opposite direction $-\mathbf{u}_k$. The step reducing strategy is performed as follows: the new step size is $\lambda_k = \varepsilon_k \lambda_{k-1}$ and the reduced factor sequence is

$$\varepsilon_k = \begin{cases} 1/1.1, & 0 < k \leq 6, \\ 1/1.2, & 6 < k \leq 8, \\ 1/1.5, & 8 < k \leq 10, \\ 1/2, & 10 < k \leq 16, \\ 1/5, & 16 < k \leq 20, \\ 1/10, & 20 < k \leq 40, \\ 1/100, & 40 < k \leq 50. \end{cases}$$

The step reducing strategy for the orthogonal shift direction must not decrease the size of shift to zero. If we cannot find a feasible point for the orthogonal shift any random feasible point can be selected.

Construction of conjugate directions with the help of orthogonal shift allows effectively searching along the boundaries of inequality constraints. Fig. 2 illustrates the search along the linear boundaries of a constraint.

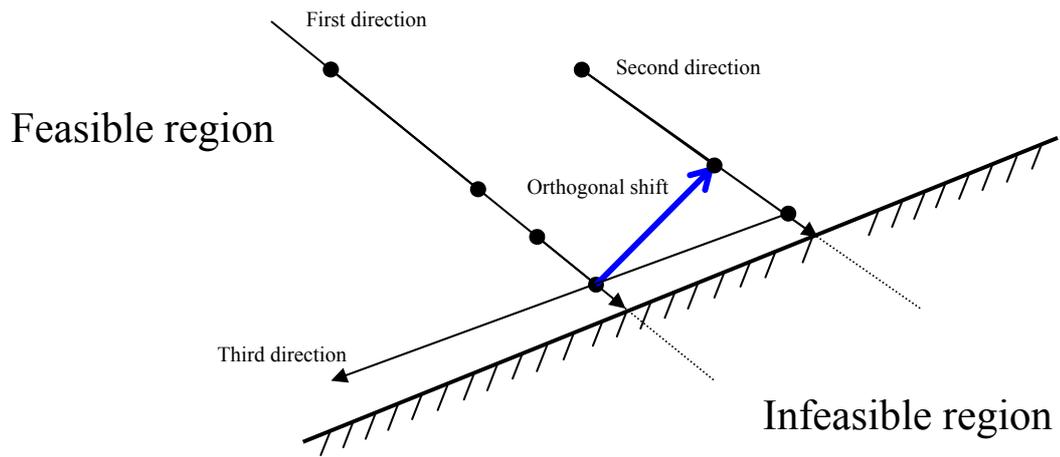

**Fig. 2. Searching along the constraint boundary**

Fig. 3 shows the real path of the following function minimum search

$$f(x, y) = x + 10y,$$

subject to

$$y \leq 2x, \ y \geq \frac{x}{2},$$

starting with initial point $[x = 100, \ y = 75]$. The real search is performed by Search command from DirectSearch package, version 2 [1]. The infeasible regions are showed in grey. The red circles denote the feasible trial points.

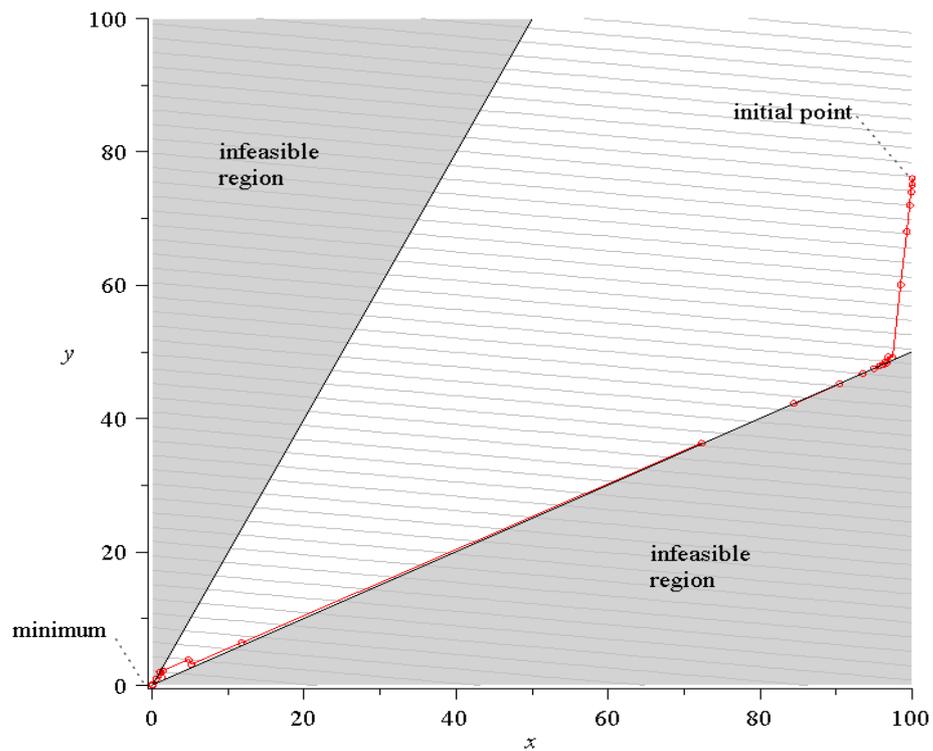

**Fig. 3. Search path**

The figure above shows the effectiveness of the CDOS method for optima searching along linear boundaries of constraints. Now we consider more complicated non-linear constraint boundaries.

Fig. 4 shows the real path of the following function minimum search

$$f(x,y) = |y-x|^{2.07} + |xy|^{1.07},$$

subject to

$$\{x \leq -1,\ x \geq -17.001,\ y \leq -1,\ y \geq -\frac{1}{3}x - 28,\ (y+20)^2 - 3x \geq 51,$$
$$|x+14.5| + (y+15)^2 \geq 3,\ (x+16)^2 + |y+8|^{3/2} \geq 20,$$
$$(x+9.2)^2 + |y+12| \geq 7,\ (x+6)^2 + (y+15)^2 \geq 29.8,\ (x+6)^2 + |y+1|^{3/2} \geq 15\},$$

starting with initial point $[x = -1.1,\ y = -27]$. This function has one minimum $f_{min} = 1$ in point $\mathbf{x}_{min} = (-1, -1)$. The constraints were chosen to make the minimum searching maximally difficult. The infeasible regions are showed in grey. The red circles denote the feasible trial points.

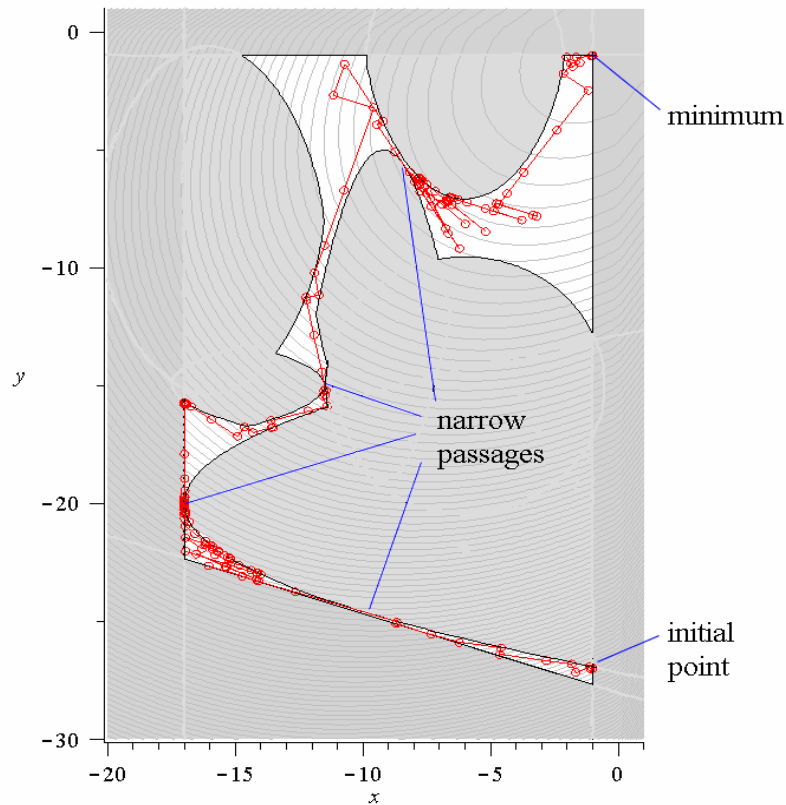

**Fig. 4. Search path**

The figure above shows the effectiveness of the CDOS method for optima searching along non-linear boundaries of constraints.

Now we consider more complicated tasks with zigzag-like feasible region.

Fig. 5 shows the real path of the following function minimum search

$$f(x,y) = \frac{|x-100|}{200} + |y-101|,$$

subject to

$$\{x \leq 100,\ x \geq 0,\ y \leq 101.01,\ y \geq 0,$$
$$|x| + |y-5|^{7/2} \geq 99.9,\ |x-100| + |y-12|^3 \geq 99.9,$$
$$|x| + |y-19|^{7/2} \geq 99.9,\ |x-100| + |y-26|^3 \geq 99.9,$$
$$|x| + |y-33|^{7/2} \geq 99.9,\ |x-100| + |y-40|^3 \geq 99.9,$$
$$|x| + |y-47|^{7/2} \geq 99.9,\ |x-100| + |y-54|^3 \geq 99.9,$$
$$|x| + |y-61|^{7/2} \geq 99.9,\ |x-100| + |y-68|^3 \geq 99.9,$$
$$|x| + |y-75|^{7/2} \geq 99.9,\ |x-100| + |y-82|^3 \geq 99.9,$$
$$|x| + |y-89|^{7/2} \geq 99.9,\ |x-100| + |y-96|^3 \geq 99.9\},$$

starting with initial point $[x=0,\ y=0]$. This function has one minimum $f_{min} = 0$ in point $x = 100,\ y = 101$. The constraints were chosen to make the minimum searching maximally difficult. The infeasible regions are showed in yellow. The red circles denote the feasible trial points.

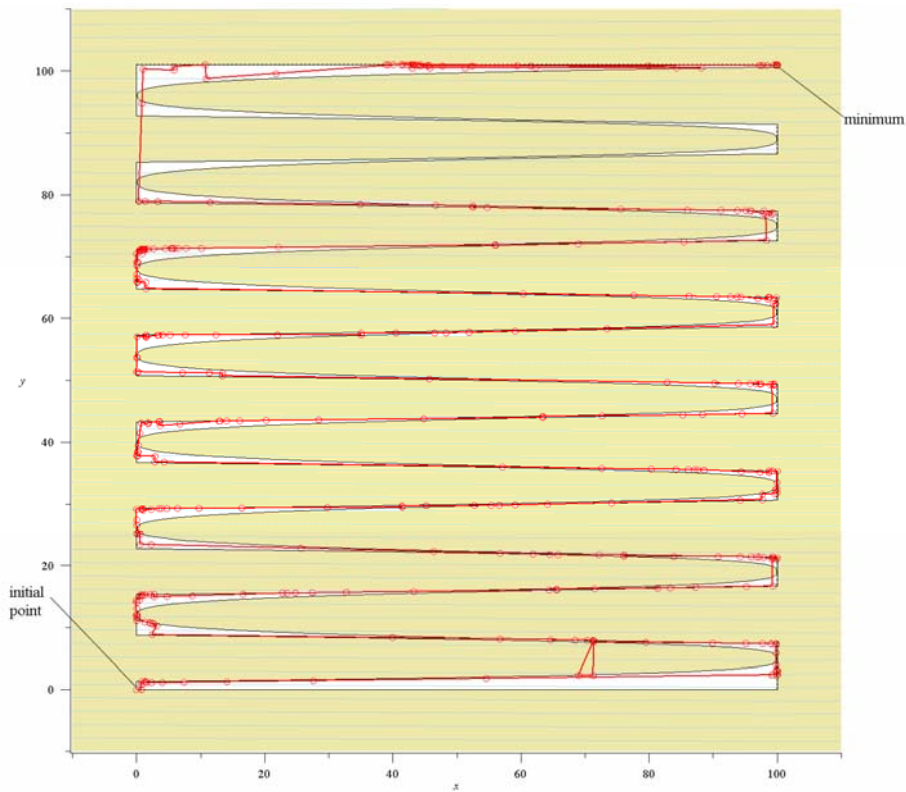

**Fig. 5. Search path**

The figure above shows the effectiveness of the CDOS method for optima searching along complex zigzag-like non-linear boundaries of constraints.

For equality constraints CDOS uses a penalty function method with a quadratic penalty function.

## Global optima search

The CDOS method can be easily transformed into an effective global optima search method because it is not a greedy searcher. For this purpose one should use the following three rules.

1) Use a large search step size, for example 100. The step size must be much greater than a size of a typical basin of attraction of the local optima.

2) The step size reducing strategy must be slower than for the local searcher. For local searching the step adaptation strategy is $\lambda^{(i+1)} = 0.3 \| \mathbf{x}_{min}^{(i)} - \mathbf{x}_{min}^{(i-1)} \| + 0.091 \cdot \lambda^{(i)}$, to search for global optima the step adaptation strategy may be $\lambda^{(i+1)} = \| \mathbf{x}_{min}^{(i)} - \mathbf{x}_{min}^{(i-1)} \| + 0.82 \cdot \lambda^{(i)}$. If you approximately know the minimum distance $d_{min}$ between two neighboring local extremum points then you can use a more effective step adaptation strategy. When the current step size is greater than $\frac{d_{min}}{2}$ the global step adaptation strategy is used, otherwise, the local step adaptation strategy is used.

3) Use a multistart strategy, i.e., use many initial points, for example 50 or more.

## Multimodal optimization

The goal of multimodal optimization is to find all local and global optima (as opposed to the single best optima).

The CDOS method can be easily transformed to an effective multimodal optimizer because it has quadratic convergence. For this purpose one should use the following three rules.

1) Use a small search step size, for example 0.005, in order not to miss the local optima. The step size must be less than the size of a typical basin of attraction of the local optima.

2) The step size reducing strategy must be the same as for the local searcher. For local searching the step adaptation strategy is $\lambda^{(i+1)} = 0.3 \| \mathbf{x}_{min}^{(i)} - \mathbf{x}_{min}^{(i-1)} \| + 0.091 \cdot \lambda^{(i)}$.

3) Use a multistart strategy, i.e., use many initial points, for example 100 or more.

Because the CDOS method has quadratic convergence it will quickly converge to the local optima as compared to methods that do not have quadratic convergence.

## Mixed integer-discrete-continuous optimization

A mixed integer-discrete-continuous optimization is a difficult task. It is partly similar to the global optimization task.

The CDOS method can be easily transformed into an effective mixed optimizer that can solve many kinds of mixed integer-discrete-continuous optimization problems. For this purpose one should use the following three rules.

1) Round the current search trial point to the nearest discrete point.
2) Use a large search step size, for example 100. The step size must be much greater than the maximum distance between two discrete neighboring points.
3) The step size reducing strategy must be slower than for the local searcher. For local searching the step adaptation strategy is $\lambda^{(i+1)} = 0.3 \| \mathbf{x}_{min}^{(i)} - \mathbf{x}_{min}^{(i-1)} \| + 0.091 \cdot \lambda^{(i)}$, for mixed optimization the step adaptation strategy may be $\lambda^{(i+1)} = \| \mathbf{x}_{min}^{(i)} - \mathbf{x}_{min}^{(i-1)} \| + 0.51 \cdot \lambda^{(i)}$. If you know the minimum distance $d_{min}$ between two neighboring discrete points (for example, $d_{min} = 1$ for integer points) then you can use a more effective step adaptation strategy. When the current step size is greater than $\frac{d_{min}}{2}$ the mixed optimization step adaptation strategy is used, otherwise, the local step adaptation strategy is used.

## Comparison with other derivative-free methods with quadratic convergence

The CDOS method was tested using the set of more than 100 test functions, for more than 500 initial points for each test function and for various method parameters and options. The method proves to be fast and reliable.

The CDOS method was compared to the following derivative-free methods with quadratic convergence:

- Powell's conjugate direction method [2, 3],
- Brent's principal axis method [4],
- Quadratic (successive quadratic approximation) [5].

The CDOS, Powell's and Brent's method are line-search optimization methods. The Quadratic method performs successive quadratic approximation of the objective function. The one-dimensional line search is performed between two successive extremum points of quadratic approximation. The quadratic approximation method is the fastest derivative-free optimization method for quadratic functions. It is similar to the Newton–Raphson method but without using the derivatives. The Quadratic method requires minimal possible number of function evaluations $\frac{(n+1)(n+2)}{2} + 1$ to exactly find the optima of the $n$-dimensional quadratic function.

The realizations of all four methods were taken from the DirectSearch package, version 2 [1].

Comparison of various methods of multivariable function optimization is difficult because of the effects of deterministic chaos during optimization. Even insignificant variation of optimization parameters may cause significant variation in the number of objective function evaluations and accuracy of the optimum value approximation. Therefore to compare it is necessary to use averaged performances of optimization efficiency.

We will characterize the optimization efficiency by the number of objective function evaluations, accuracy of objective function optimum evaluation and reliability. The reliability is described by the percentage of the number of successful optimizations.

For all methods we take the same key optimization parameters: the step size is equal to 1, the tolerances for optima and optimum point are $10^{-6}$. When the absolute difference between the optimum value and its approximation is greater than "fmin maximum"=$10^{-3}$ we decide that the optimization has failed.

## *Differentiable functions*

For quadratic and near quadratic functions all four methods prove to be fast and reliable because all methods have quadratic convergence for such functions.

For differentiable but non-quadratic functions the line-search methods as a rule outperform the Quadratic method. To show typical performances of all four methods for differentiable but non-quadratic functions we select the Rosenbrock test function to minimize

$$f(x, y) = 100(y - x^2)^2 + (1 - x)^2.$$

This function has one minimum $f_{\min} = 0$ in point $x = 1, y = 1$. We take the following 500 initial points:

$$(x_i = -1 + i,\ y_i = 2 + i),\ i = 0, 1, 2, ..., 499.$$

Fig. 6 shows the number of function evaluations by four methods. Fig. 7 shows the accuracy of objective function optimum evaluation.

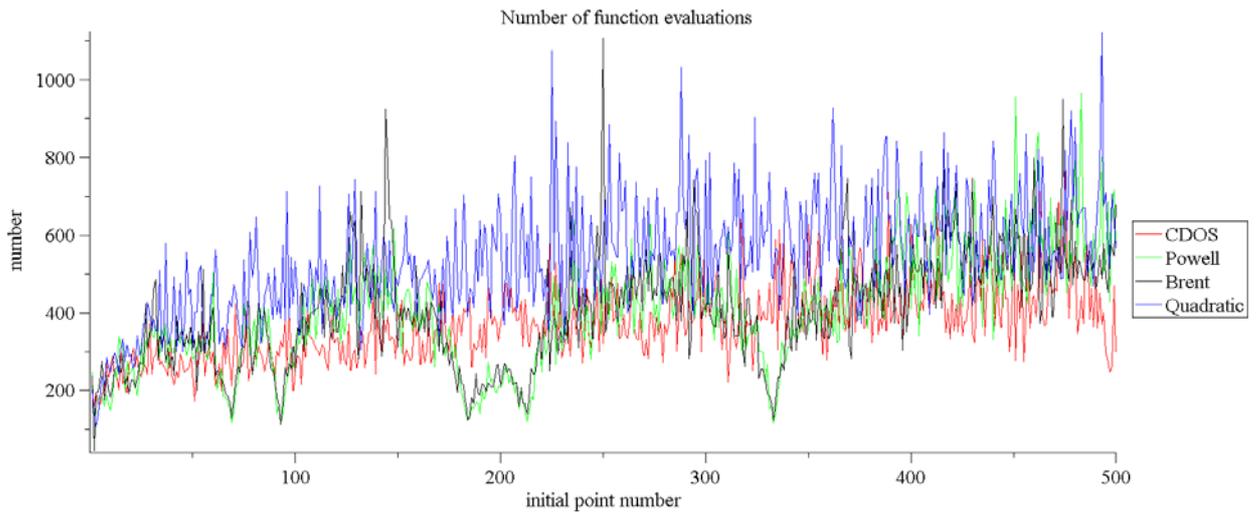

**Fig. 6. Number of function evaluations**

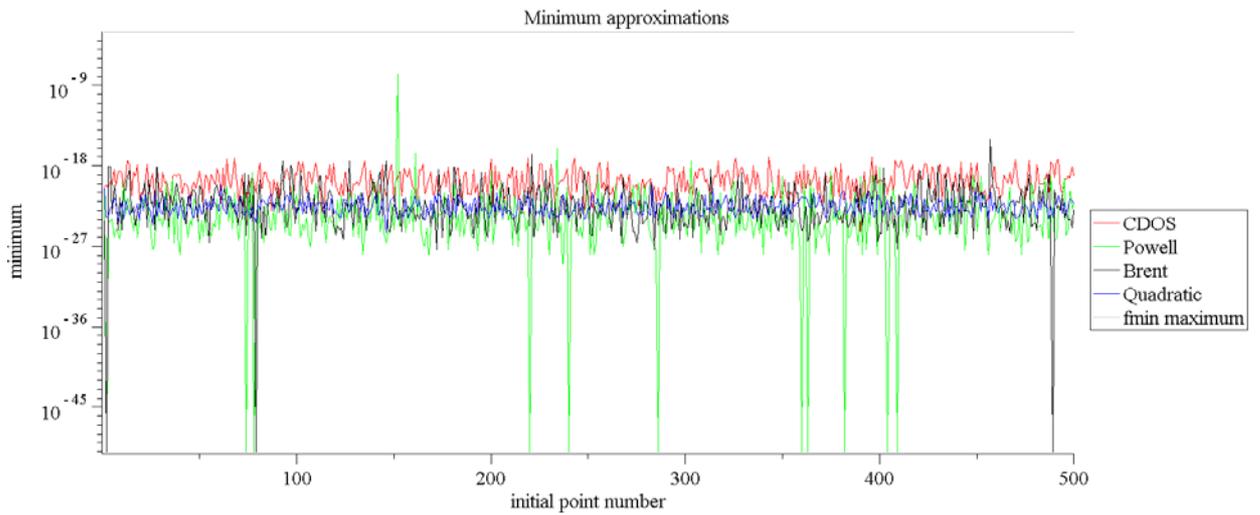

**Fig. 7. Minimum function approximation**

Table 1 shows the total performances of four methods for the Rosenbrock test function.

**Table 1. Total efficiency of optimization methods**

| Method | CDOS | Powell | Brent | Quadratic |
|---|---|---|---|---|
| Average number of function evaluations | 372 | 402 | 403 | 527 |
| Median of minimum approximation | $10^{-20}$ | $10^{-24}$ | $10^{-23}$ | $10^{-23}$ |
| Reliability | 100% | 100% | 100% | 100% |

Fig. 7 and Table 1 show that the quadratic convergence of four methods leads to a much better accuracy of minimum approximation than the tolerance given.

The optimization effectiveness of the line-search methods is approximately the same. The Powell's and Brent's methods may be faster for some very bad scaled functions, and the CDOS method is more reliable and faster when the initial point is far away from the optima point.

## Non-differentiable functions

The CDOS method is reliable for non-differentiable functions, while Powell's, Brent's and Quadratic methods turn out to be very unreliable for such functions. The following example demonstrates this fact.

We select the non-differentiable function similar to the Rosenbrock test function to minimize

$$f(x,y) = 100|y - x^2| + |1 - x|.$$

This function has one minimum $f_{min} = 0$ in point $x = 1$, $y = 1$. We take the following 500 initial points:

$$(x_i = -1 + i,\ y_i = 2 + i),\ i = 0, 1, 2, ..., 499.$$

Fig. 8 shows the number of function evaluations by four methods. Fig. 9 shows the accuracy of objective function optimum evaluation.

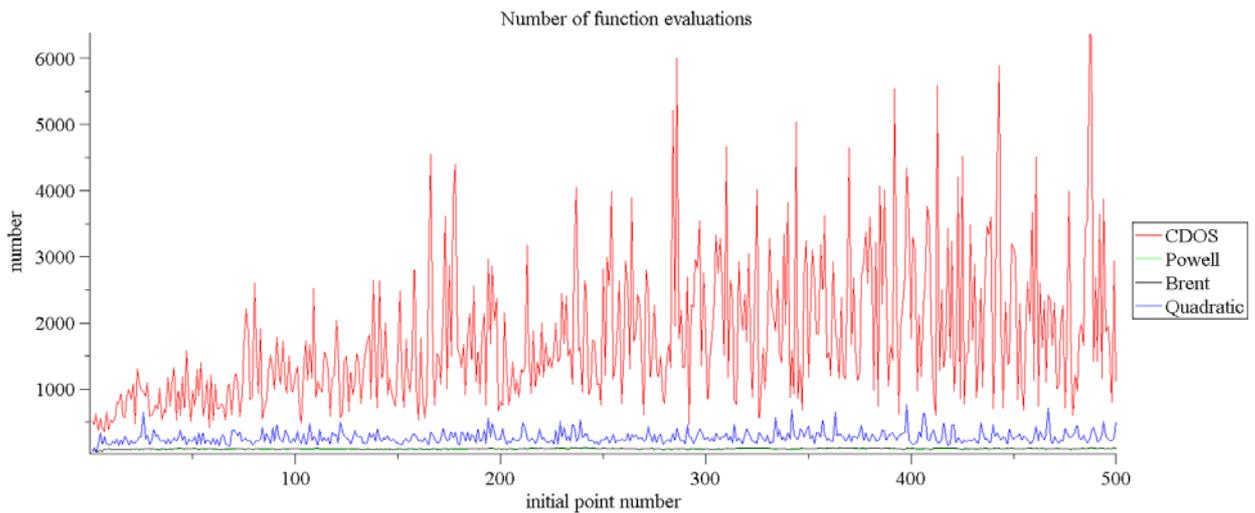

**Fig. 8. Number of function evaluations**

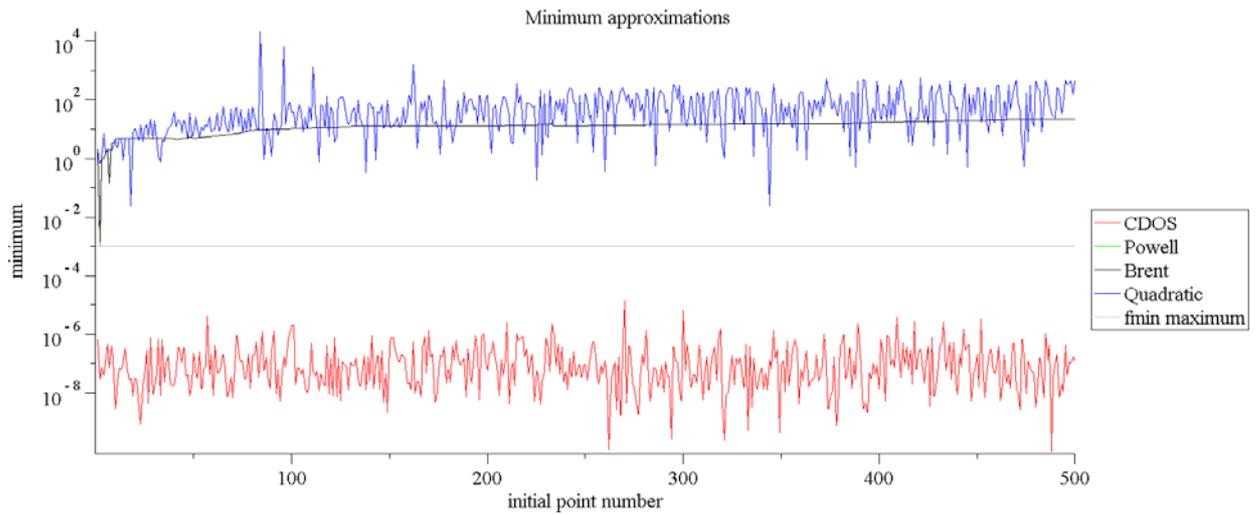

**Fig. 9. Minimum function approximation**

Table 2 shows the total performances of four methods for the non-differentiable function similar to the Rosenbrock test function.

**Table 2. Total efficiency of optimization methods**

| Method | CDOS | Powell | Brent | Quadratic |
|---|---|---|---|---|
| Average number of function evaluations | 1751 | 99 | 102 | 274 |
| Median of minimum approximation | $7 \cdot 10^{-8}$ | 13.5 | 13.5 | 30.5 |
| Reliability | 100% | 0% | 0% | 0% |

Fig. 9 and Table 2 show that the reliability of Powell's, Brent's and Quadratic methods is zero for this function.

It is known that the Nelder-Mead simplex method is good enough for non-differentiable functions. This method copes with the non-differentiable function similar to the Rosenbrock test function but it requires approximately 15 times more function evaluations than the CDOS method.

## *Constrained optimization*

The CDOS method is reliable to solve the constrained optimization tasks when inequality constraints are handled directly while the Powell's, Brent's and Quadratic methods turn out to be unreliable for such tasks. The following example proves this.

We select the following simple function to minimize

$f(x, y) = x + 10y$,

subject to

$$y \leq 2x, \ y \geq \frac{x}{2}.$$

This function has one minimum $f_{min} = 0$ in point $x = 0, \ y = 0$. We take the following 500 initial points:

$$(x_i = i, \ y_i = i), \ i = 1, 2, ..., 500.$$

Fig. 10 shows the number of function evaluations by four methods. Fig. 11 shows the accuracy of objective function optimum evaluation.

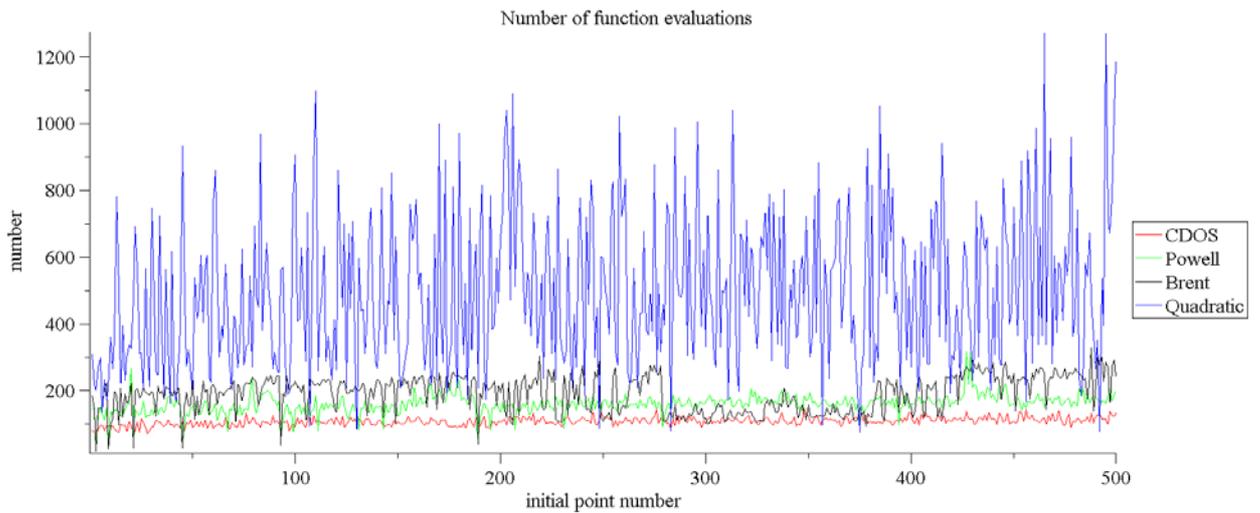

**Fig. 10. Number of function evaluations**

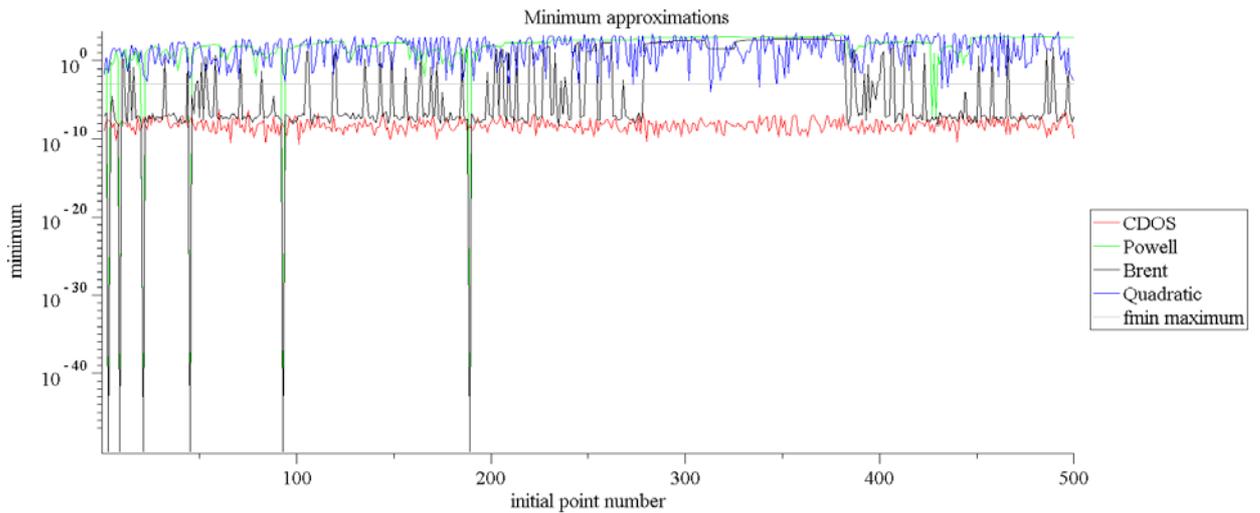

**Fig. 11. Minimum function approximation**

Table 3 shows the total performances of four methods for the constrained optimization.

Table 3. Total efficiency of optimization methods

| Method | CDOS | Powell | Brent | Quadratic |
|---|---|---|---|---|
| Average number of function evaluations | 108 | 160 | 191 | 493 |
| Median of minimum approximation | $6 \cdot 10^{-9}$ | 202 | $10^{-7}$ | 48.8 |
| Reliability | 100% | 1.8% | 62.6% | 1.2% |

Table 3 shows that the reliability of Powell's and Quadratic methods is very low for this function. The reliability of Brent's method is greater but insufficient for such simple optimization task. The CDOS method shows high reliability for constrained optimization problems.

## Conclusion

It is known that there are no universal optimization methods which can effectively solve different kinds of optimization problems. We call the CDOS method universal because it proves to be fast and reliable for a wide range of optimization problems. It is fast and reliable not only for differentiable functions but also for non-differentiable functions as well as constrained optimization tasks where the objective function may be undefined outside the constraint boundaries and both the objective function and constraints may be black boxes.

## References

1. S.N. Moiseev, "DirectSearch optimization package, version 2", Maple application centre, http://www.maplesoft.com/applications/view.aspx?SID=101333. Feb 1, 2011.

2. Powell, M. J. D., 1964, An efficient method for finding the minimum of a function of several variables without calculating derivatives, *Computer J.*, 7,155-162.

3. Himmelblau M. David. Applied Nonlinear Programming. McGraw-Hill Book Company. 1972. 534 p.

4. Richard P. Brent, Algorithms for Minimization without Derivatives, Prentice-Hall, Englewood Cliffs, NJ, 1973, 195 p.

5. Numerical methods for constrained optimization. Chapter 7. / Edited by P.E. Gill and W. Murrey, London; New York: Academic Press, 1974, 283 p.